\documentclass{article}

\usepackage{natbib}
\usepackage{amssymb}
\usepackage{amsfonts}
\usepackage{latexsym}
\usepackage[latin1]{inputenc}
\usepackage[english]{babel}

\begin{document}


\title{Computing Bayesian predictive distributions:
The K-square and K-prime distributions}


\maketitle

\begin{center}
{\bf Jacques Poitevineau$^{1}$, \bf Bruno Lecoutre$^2$}\\
\end{center}

\noindent $^1$ERIS and LAM/LCPE, UMR~7604, C.N.R.S.,
Université Paris~6 et Ministère de la Culture,
11 rue de Lourmel, 75015~Paris, France.\\
Email: poitevin@ccr.jussieu.fr\\

\vspace*{-0.2truecm}
\noindent $^2$ERIS, and Laboratoire de Mathématiques Raphaël Salem,\\
UMR~6085, C.N.R.S. et Université de Rouen, Avenue de l'Université, \\
BP 12, 76801~Saint-Etienne-du-Rouvray, France.\\
E-mail: bruno.lecoutre@univ-rouen.fr\\

\begin{abstract}
The computation of two Bayesian predictive distributions which are 
discrete mixtures of incomplete beta functions is considered. The
number of iterations can easily become large for these distributions 
and thus, the accuracy of the result can be questionable. Therefore, 
existing algorithms for that class of mixtures are improved by 
introducing round-off error calculation into the stopping rule. 
A further simple modification is proposed to deal with possible 
underflows that may prevent recurrence to work properly.
\end{abstract}



\section{Introduction}

The K-square and K-prime distributions have been introduced in
\citet{BL84}. They can be characterized as mixtures of
the classical noncentral $F$ and noncentral $t$ distributions 
respectively \citep{BL99}. These two distributions are involved in 
the Bayesian predictive approach for planning and monitoring 
experiments \citep{BL01}. In particular, they are useful tools for 
sample size determination, using the predictive distributions of the
test statistics and of the limits of confidence intervals under 
standard normal models, assuming a conjugate prior. It must also be 
noted that they include as particular cases the distributions of the 
square of the sample multiple correlation coefficient and of the 
sample correlation coefficient. The aim of this article is to 
provide efficient algorithms for the calculation of their cumulative 
distribution functions (cdfs). These cdfs can be expressed in terms 
of infinite series of multiples of incomplete beta function ratios, 
thus adequate for recursive calculations. More precisely, both imply 
the general form
\begin{equation}
\sum_{j=0}^\infty s^jg_jH_j(x),
\label{EQ1}
\end{equation}
with 
\[ s=\pm 1,\quad 0\leq g_j\leq 1\ \forall j,
\quad \sum_{j=0}^\infty g_j=1\]
and where $H_j(x)$ involves only the incomplete beta function.

Dealing with a related problem, the Applied Statistics algorithm
AS 278 developed for the psi-square distribution \citep{psi2} could
be adapted to match the present cdfs. However, AS 278 is a Method~1
recursive algorithm, in the terms of \citet{BK03}: accumulation is
simply done from index $0$ until a convergence criterion is met. In
some cases (especially when the noncentrality parameter of
the distribution is large), it can lead to an exceedingly large number
of iterations, and consequently to unacceptable execution time and
loss of precision. \citet{Frick} proposed an improvement that consists
in starting iterations at an index such that the resulting truncation
error is negligible, but this does not solve the problem.

Yet, the present cdfs are of the general class considered
by \citet{BK03} and, as such, are good candidates for what they called
Method~2 class of algorithms. Essentially, this Method~2 is a
both backward and forward recursive algorithm where the starting index
for iterations, say $k$, is chosen so that $g_k$ is a maximum, which
reduces the above mentioned problems. Nevertheless, although smaller 
than with Method~1, the number of iterations can  still remain 
important as soon as parameters increase. Thus, when a relatively high
degree of accuracy is required, the problem of round-off errors cannot
be neglected.

Therefore, we present in the next two sections a Method~2 class of
algorithms that includes round-off error calculations. It is applied
here respectively to the K-square and K-prime cdfs, but is of general
use as far as the general form (\ref{EQ1}) is concerned. 
CPU times are presented in section \ref{secCPU}, along with 
a few illustrations, and some examples of applications of these cdfs 
are given in section \ref{secEX}. In section \ref{secPBI}
we discuss some remaining problems and propose, in some cases, a
simple modification which leads to an algorithm that is intermediate 
between Method~1 and Method~2. Section \ref{secCLR} is devoted to some
concluding remarks.

\section{K-square distribution} \label{secK2}

Technical characterizations of the K-square distribution can be found 
in \citet{BL99}. This distribution is written $K^2_{p,q,r}(a^2)$
where $p, q, r$ are degrees of freedom parameters and $a^2$ is a
noncentrality parameter.

Particular cases of the K-square distribution are:\\
$a=0$ : \ $K^2_{p,q,r}(0) \equiv F_{p,r}$ (usual $F$ distribution),\\
$q=\infty$ : \ $K^2_{p,\infty,r}(a^2) \equiv F'_{p,r}(a^2)$
(noncentral $F$ distribution),\\
$r=\infty$ : \ $K^2_{p,q,\infty}(a^2) \equiv \Lambda_{p,q}^2(a^2)$
(lambda-square or alternate chi-square distribution),\\
$q=\infty, r=\infty$ : \ $K^2_{p,\infty,\infty}(a^2)
\equiv (1/p)\chi_p^2(a^2)$ (noncentral chi-square distribution).\\

For the cdf, $s=1$ in (\ref{EQ1}) and we simply have
\[ \Pr(K_{p,q,r}^2(a^2)<x)=\sum_{j=0}^{\infty} g_jH_j(x), \]
with
\begin{equation}
g_j=\frac{\Gamma(\frac{q}{2}+j)}{\Gamma(j+1)\Gamma(\frac{q}{2})}
\left(\frac{q}{q+a^2}\right)^{\frac{q}{2}}
\left(\frac{a^2}{q+a^2}\right)^j
\end{equation}
and
\begin{equation}
H_j(x)=I_{px/(r+px)}\left(\frac{p}{2}+j, \frac{r}{2}\right),\ x>0,
\end{equation}
where $I_z$ is the incomplete beta function
\[
I_z(a,b)=\frac{\Gamma(a+b)}{\Gamma(a)\Gamma(b)}
\int_{0}^{z} t^{a-1}(1-t)^{b-1}dt.
\]

The coefficients $g_j$ are the probabilities of obtaining the value
$j$ for a variate following a negative binomial distribution with
parameters $q/(q+a^2)$ and $q/2$. The mode is $[a^2(q-2)/(2q)]$,
where $[.]$ denotes the integer part \citep[e.g., see][p. 209]{JKK93},
hence the starting index for iterations. From this, it is clear that
the number of iterations heavily depends on $a^2$.

The recurrence relations for the cdf are straightforward. For the
$H_j$'s (the incomplete beta function) we have
\begin{eqnarray*}
H_{j+1}&=&H_j-\frac{\Gamma(p/2+r/2+j)}{\Gamma(p/2+j+1)\Gamma(r/2)}
\ \left(\frac{px}{r+px}\right)^{p/2+j}
\left(\frac{r}{r+px}\right)^{r/2},\\
H_{j-1}&=&H_j+\frac{\Gamma(p/2+r/2+j-1)}{\Gamma(p/2+j)\Gamma(r/2)}
\ \left(\frac{px}{r+px}\right)^{p/2+j-1}
\left(\frac{r}{r+px}\right)^{r/2}
\end{eqnarray*}
and for the $g_j$ coefficients
\begin{eqnarray*}
g_{j+1}&=&\frac{q/2+j}{j+1}\ \frac{a^2}{q+a^2}\quad g_j,\\
g_{j-1}&=&\frac{j}{q/2+j-1}\ \frac{q+a^2}{a^2}\quad g_j.
\end{eqnarray*}

Let $\Delta$ and $\delta$ denote the absolute and the relative 
error respectively. The absolute error for an individual term of the 
series is
\[ \Delta (g_jH_j)=g_j\Delta H_j\ +\ H_j\Delta g_j. \]
Now, noting $k$ the starting index of the computations, the forward
and backward recurrences for $g_j$ are respectively of the form
\[ g_{k+j}=g_{k+j-1}c_{k+j-1}=g_k \prod_{i=0}^{j-1}c_{k+i}
\quad \textrm{and} \quad g_{k-j}=g_{k-j+1}/c_{k-j}=g_k
\prod_{i=1}^{j}\frac{1}{c_{k-i}}, \]
so that
\[ \delta	g_{k+j}=\delta g_k +  \sum_{i=0}^{j-1}\delta c_{k+i}
\quad \textrm{and} \quad
   \delta	g_{k-j}=\delta g_k +  \sum_{i=1}^{j}\delta c_{k-i}.\]
If we assume that the relative errors on the coefficients $c_j$
are constant, say equal to $\epsilon$ (e.g., we can assume that all
$c_j$'s are calculated with a maximal precision of $n$ decimal digits
so that $\epsilon < \frac{1}{2}10^{-n+1}$), we obtain
\[ \delta	g_{k\pm j}=\delta g_k + j\epsilon \quad \textrm{hence}\quad
 \Delta g_{k\pm j}=(\delta g_k + j\epsilon)g_{k\pm j}. \]
For the terms $H_j(x)$, the recurrence involves a sum
\[ H_{k+j}=H_{k+j-1}-d_{k+j-1}=H_k-\sum_{i=0}^{j-1}d_{k+i}, \]
\[ H_{k-j}=H_{k-j+1}+d_{k-j}=H_k+\sum_{i=1}^jd_{k-i}, \]
then,
\[ \Delta H_{k+j}=\Delta H_k\ +\ \sum_{i=0}^{j-1}\Delta d_{k+i}
\ \quad \textrm{and} \quad
\Delta H_{k-j}=\Delta H_k\ +\ \sum_{i=1}^j\Delta d_{k-i}. \]
The coefficients $d_{k\pm i}$ contain gamma functions which can
themselves be calculated by recurrence, just as
for the $g_j$'s. Therefore, with the same assumptions as for the
coefficients $g_j$, we have
\[ \Delta d_{k\pm j}=(\delta d_k + j\epsilon)d_{k\pm j}. \]
Consequently, the round-off error ($E_c$) of a calculation involving
$N$ iterations (both backward and forward) becomes
\begin{eqnarray}
E_c&=&\Delta (g_kH_k)+\sum_{j=1}^N\Delta(g_{k+j}H_{k+j})+
\sum_{j=1}^{\min(N,k)}\Delta(g_{k-j}H_{k-j}).
\end{eqnarray}
Skipping tedious but elementary calculations, it gives
\begin{eqnarray}
E_c&=&(\delta g_k+\delta H_k)g_kH_k+ \nonumber \\
   &&  \sum_{j=1}^N\{(\delta H_k+\delta d_k)g_{k+j}H_k+
   (\delta g_k+j\epsilon-\delta d_k)g_{k+j}H_{k+j}+ \nonumber \\
   && \epsilon g_{k+j}\sum_{i=0}^{j-1}id_{k+i}\}+ \nonumber \\
   &&  \sum_{j=1}^{\min(N,k)}\{(\delta H_k-\delta d_k)g_{k-j}H_k+
   (\delta g_k+j\epsilon+\delta d_k)g_{k-j}H_{k-j}+ \nonumber \\
   &&  \epsilon g_{k-j}\sum_{i=1}^{j}id_{k-i}\}.
\end{eqnarray}
Now, for the same reason as for the relative errors on the 
coefficients $c_j$, we can assume $\delta g_k=\delta d_k=\epsilon$. 
Furthermore, $H_k$ involves only one calculation of the incomplete 
beta function for which there exist very performing algorithms 
\citep[e.g., AS 63 by][]{AS63}, so that, again, 
$\delta H_k=\epsilon$ is a reasonable assumption. Consequently, 
it reduces finally to
\begin{eqnarray}
E_c&=&\epsilon \left[
	 2H_k\sum_{j=0}^Ng_{k+j}+\sum_{j=1}^Njg_{k+j}H_{k+j}+
	 \sum_{j=1}^N\left\{g_{k+j}\sum_{i=0}^{j-1}id_{k+i}\right\}
	 + \right. \nonumber \\
	 && \left. 2\sum_{j=1}^{\min(N,k)}g_{k-j}H_{k-j}+
   \sum_{j=1}^{\min(N,k)}jg_{k-j}H_{k-j}+ \right. \nonumber \\
   && \left. \sum_{j=1}^{\min(N,k)}\left\{g_{k-j}
   \sum_{i=1}^{j}id_{k-i}\right\}
   \right].
\end{eqnarray}

Given the parameters, $H_j(x)$ is a decreasing function of $j$.
Thus, when stopping the calculations at step $j$, the truncation
error ($E_t$) is bounded by: \\
\\
while $j < k$
\begin{eqnarray}
E_t &\leq& H_0(x)\sum_{i=0}^{k-j-1}g_i +
			H_k(x)\sum_{i=k+j+1}^\infty g_i \nonumber \\
    &\leq& H_0(x)\sum_{i=0}^{k-j-1}g_i +
    	H_0(x)\sum_{i=k+j+1}^\infty g_i \nonumber \\
    &\leq& H_0(x)\left[ 1-\sum_{i=k-j}^{k+j}g_i\right] \label{truncH0}
\end{eqnarray}
and when $j \geq k$
\begin{equation}
E_t \leq H_{k+j}(x)\left[ 1-\sum_{i=0}^{k+j}g_i\right].
\end{equation}

(\ref{truncH0}) is a slight modification of the rule in step 3 in
\citet{BK03} who used $1$ instead of $H_0(x)$. The relaxation of the
stopping rule compensates for the increased execution
time due to one call to the incomplete beta function.

\emph{Stopping rule}: Stop when $E_t+E_c$ becomes lower than a
predetermined absolute error bound or when $E_c$ exceeds that error
bound, which means that the required accuracy cannot be reached.

For the distribution of the square of the sample multiple 
correlation coefficient (see end of section \ref{secEX}), we
compared the algorithm for the K-square cdf, called K2CDF, to 
\citet{BK03} Algorithm 7.1 (the mode of the negative binomial 
distribution, instead of the mean, was used as the starting point 
to ensure the comparability of the two algorithms). For the 
examples in their Table 1, all results agreed within the 
$10^{-12}$ limit that was chosen as the maximum absolute error 
parameter (both algorithms were run in ``double precision'',
i.e. 64-bit words).

\section{K-prime distribution} \label{secKP}

Technical characterizations of the K-prime distribution can be found
in \citet{BL99}. This distribution is written $K'_{q,r}(a)$ where
$q, r$ are degrees of freedom parameters and $a$ is a noncentrality
parameter.

Particular cases of the K-prime distributions are:\\
$a=0$ : \ $K'_{q,r}(0) \equiv t_r$ (usual $t$ distribution),\\
$q=\infty$ : \ $K'_{\infty,r}(a) \equiv t'_r(a)$ (noncentral $t$
distribution),\\
$r=\infty$ : \ $K'_{q,\infty}(a) \equiv \Lambda'_q(a)$
(lambda-prime distribution),\\
$q=\infty,r=\infty$ : \ $K'_{\infty,\infty}(a^2) \equiv N(a,1)$
(normal distribution).\\

This cdf has the following properties:
\[\Pr(K'_{q,r}(a)<x)= \Pr(K'_{r,q}(x)>a), \]
\[\Pr(K'_{q,r}(-a)<-x)= \Pr(K'_{q,r}(x)>a),\]
\[\Pr(K'_{q,r}(a)<0)= \Pr(\Lambda'_q(a)<0)=\Pr(t_q>a).\]

Several cases are to be distinguished for the cdf:\\
\\
If $a>0$ and $x<0$
\begin{eqnarray*}
\Pr(K'_{q,r}(a)<x)&=& \Pr(K'_{q,r}(a)<0)-\Pr(x<K'_{q,r}(a)<0)\\
   &=& \Pr(t_q>a)-\sum_{j=0}^{\infty} (-1)^jg_jI_{x^2/(r+x^2)}
       \left(\frac{j+1}{2}, \frac{r}{2}\right),\
\end{eqnarray*}
where
\begin{equation}
g_j=\frac{1}{2}\frac{\Gamma(\frac{q+j}{2})}{\Gamma(\frac{1+j}{2})
\Gamma(\frac{q}{2})}\left(\frac{q}{q+a^2}\right)^{\frac{q}{2}}
\left(\frac{a^2}{q+a^2}\right)^{\frac{j}{2}}.
\end{equation}
If $a>0$ and $x>0$
\begin{eqnarray*}	
\Pr(K'_{q,r}(a)<x)&=& \Pr(K'_{q,r}(a)<0)+\Pr(0<K'_{q,r}(a)<x)\\
   &=& \Pr(t_q>a)+\sum_{j=0}^{\infty} g_jI_{x^2/(r+x^2)}
   \left(\frac{j+1}{2}, \frac{r}{2}\right).
\end{eqnarray*}
If $a<0$, we reduce to the above cases using
\[\Pr(K'_{q,r}(a)<x) = 1-\Pr(K'_{q,r}(-a)<-x).\]
If $a=0$, we simply have
\[\Pr(K'_{q,r}(0)<x) = \Pr(t_r<x).\]

Hence, the cdf of the K-prime involves the calculation of the cdf of
the usual Student's $t$ distribution and a series of the general form
(\ref{EQ1}). The case where $a$ and $x$ are of a different sign is an
unfavorable one, since the series is then alternate. Therefore, in the
algorithm called KPRIMECDF, the even and odd terms of the series are 
accumulated separately in order to minimize the number of subtractions.

The recurrence relations for the incomplete beta function now write
\begin{eqnarray*}
H_{j+2}&=&H_j-\frac{\Gamma(\frac{j+r+1}{2})}
{\Gamma(\frac{j+3}{2})\Gamma(\frac{r}{2})}
\ \left( \frac{x^2}{r+x^2}\right)^{\frac{j+1}{2}}
\left(\frac{r}{r+x^2}\right)^\frac{r}{2},\\
H_{j-2}&=&H_j+\frac{\Gamma(\frac{j+r-1}{2})}
{\Gamma(\frac{j+1}{2})\Gamma(\frac{r}{2})}
\ \left( \frac{x^2}{r+x^2}\right)^{\frac{j-1}{2}}
\left(\frac{r}{r+x^2}\right)^\frac{r}{2}
\end{eqnarray*}
and for the $g_j$ coefficients
\begin{eqnarray*}
g_{j+2}&=&\frac{q+j}{j+2}\ \frac{a^2}{q+a^2}\ g_j, \\
g_{j-2}&=&\frac{j}{q+j-2}\ \frac{q+a^2}{a^2}\ g_j.
\end{eqnarray*}

The starting point for iterations is taken as the mode of the $g_j$'s,
i.e. $k=[a^2(q-2)/q]$. Again, $a^2$ is an important factor regarding
the number of iterations. The calculation of errors developed for the
K-square series directly applies here, and the stopping rule is the
same.

\section{Numerical examples and CPU time} \label{secCPU}

Some numerical examples, also illustrating the speed of the
algorithms, are presented in Tables 1 and 2. The probabilities
presented have been calculated with a required accuracy of $10^{-4}$.
In order to estimate the loss of speed due to the calculation of
round-off errors, we also computed the cdf using only the truncation
error in the stopping rule to serve as reference CPU times. In the
last column of the tables, the time increase is expressed as a
percentage of these reference CPU times. The programs were compiled 
with the GNU g95 Fortran compiler (GCC 4.0.3, Apr. 19 2006), 
using ``standard real'' data type (i.e., 32-bit words), 
and CPU time was computed through the Fortran CPU\_TIME subroutine. 
The programs were run on an Intel M750 1.86 GHz PC 
(each calculation was computed 20,000 times in order to provide a
substantial CPU time).

On the one hand, and as easily predictable from the algorithm, 
it appears that calculation of round-off errors is time consuming. 
On the other hand, examples of its usefulness can be given. 
For that purpose, keeping the required accuracy to $10^{-4}$, we 
consider that the same algorithm run in ``double precision'' 
(64-bit words) with an accuracy parameter set to $10^{-9}$ 
provides the ``exact'' value. 
The absolute difference between this reference value and the value 
returned by the algorithm without round-off error calculations is 
termed ``error'' in the following (the ``exact'' value is reported 
in square brackets). In all these cases the algorithm with round-off 
error calculations rightly returns an error message indicating the 
required accuracy cannot be met.

For the K-square cdf:\\
$x=90,\ p=10,\ q=15,\ r=20,\ a^2=10^3 :\ error=1.7\times10^{-4}
\ ~~[0.4168]$,
$x=15,\ p=10,\ q=20,\ r=10^5,\ a^2=80 :\ error=6.0\times10^{-4}
\ ~~[0.9577]$,
$x=9,\ p=10,\ q=100,\ r=10^5,\ a^2=80 :\ error=1.2\times10^{-2}
\ ~~[0.5259]$.

For the K-prime cdf:\\
$x=100,\ q=10,\ r=20,\ a=80 :\ error=9.0\times10^{-4}\ ~~[0.8101]$, \\
$x=20,\ q=10,\ r=10^5,\ a=20 :\ error=4.9\times10^{-3}\ ~~[0.5574],$ \\
$x=20.5,\ q=200,\ r=10^6,\ a=21 :\ error=1.5\times10^{-1}\ ~~[0.3730].$

All these examples involve the largeness of at least one parameter, 
precisely because it is in such cases that the precision of the result 
may be suspected. An illustrated example for the K-prime cdf is 
presented in the next section.

\begin{table}
\caption{Time comparison between Algorithm K2CDF and the same
algorithm without round-off error calculation for computing
$\Pr(K_{p,q,r}^2(a^2)<x)$ 20,000 times (time in second)}
\begin{tabular}{l l l l l l l l}
 \hline
 $x$&$p$&$q$&$r$&$a^2$&$\Pr(K_{p,q,r}^2(a^2)<x)$&CPU time&time increase\\
 \hline
   3& 5& 5& 5&   5&0.6664& 0.20& 08\% \\
   1& 5& 5& 9&  10&0.1195& 0.11& 17\% \\
  10& 5& 5& 9&  10&0.9440& 0.14& 29\% \\
  10& 5& 5& 9& 100&0.2142& 0.25& 14\% \\
 100& 9& 5& 5& 100&0.9819& 0.53& 31\% \\
  80&10&20&25&1000&0.3015& 1.31& 27\% \\
 \hline
\end{tabular}
\end{table}

\begin{table}
\caption{Time comparison between Algorithm KPRIMECDF and the same
algorithm without round-off error calculation for computing
$\Pr(K_{q,r}'(a)<x)$ 20,000 times (time in second)}
\begin{tabular}{l l l l l l l}
 \hline
 $x$&$q$&$r$&$a$&$\Pr(K_{q,r}'(a)<x)$&CPU time&time increase\\
 \hline
  -5&5&    5&  0.5&0.0007&0.50& 07\% \\
   5&5&    5&  5  &0.5000&0.34& 10\% \\
   9&5&    5&  5  &0.8763&0.55& 25\% \\
   5&5&    5& 10  &0.0872&0.47& 15\% \\
   9&5&    5& 10  &0.4137&0.77& 26\% \\
   9&5&10000&  5  &0.9856&0.45& 16\% \\
 -15&5&   10&-50  &0.9918&4.47& 30\% \\  
 \hline
\end{tabular}
\end{table}

\section{Examples of applications} \label{secEX}
As an illustration of the use of the K-prime and K-square
distributions, consider the sample size determination under usual
normal models. For instance, a simple two-sample experiment is
designed to compare a new drug with a placebo. The goals of the
experiment specify that the new drug is considered as effective if
the raw difference $\delta=\mu_D-\mu_P$ is more than +3.
For this purpose, the investigators plan to use a two-sample shifted
$t$ test with equal numbers of subjects $n$ in each group, in order
to test H$_0: \delta=+3$ against the alternative H$_1: \delta>+3$.
Hence, the efficacy of the drug will be assessed if

\[ t = \frac{d-3}{s\sqrt{2/n}} > t_{q,0.05}, \]

\noindent where $d$ is the observed difference, $s$ is the pooled
estimate of the common standard deviation $\sigma$ and $t_{q,0.05}$
is the 5\% upper point of the Student's distribution with $q=2n-2$
degrees of freedom.

Suppose that a conjugate prior distribution has been chosen, such
as $\delta|\sigma \sim N(d_0,(2/n_0)\sigma^2)$ and $\sigma^2
\sim s_0^2(\chi_{q_0}^2)^{-1}$. For instance, this prior can be the
posterior distribution from a pilot study (starting with a
noninformative prior). Then, for any given sample size $n$, the
probability of achieving the study target can be computed from a
K-prime distribution, using the predictive distribution of the $t$
test statistic:

\[ t \sim \sqrt{1+n/n_0}\, K'_{q_0,q}\left(
\frac{t_0}{\sqrt{1+n_0/n}}\right), ~~\hbox{where }
t_0 = \frac{d_0-3}{s_0\sqrt{2/n_0}}. \]

Suppose that $d_0=+4.35$, $s_0=2.07$, $n_0=10$, hence $q_0=18$ and 
$t_0 = +1.458$. For instance we find for $n=50$ the predictive
probability:

\[ Pr(t> +1.6606) = Pr \left[ K'_{18,98}\left(1.458\sqrt{5/6}\,\right)
 > 1.6606/\sqrt 6\, \right] = 0.7327. \]

In order to get predictive probabilities equal to 0.80 and to 0.90,
$n=97$ and $n=1930$ subjects in each group are respectively needed.

Equivalently, the investigators could compute a 90\% confidence
interval for $\delta$ and assess the efficacy of the drug if its lower
limit is larger than +3. The predictive distribution for this lower
limit $\underline\ell= d - t_{q,0.05}\,s\sqrt{2/n}$ also involves a
K-prime distribution:

\[ \underline\ell \sim d_0 - s_0 \sqrt{2/n_0+2/n}
	\ K'_{q,q_0}\left(\frac{-s_0t_{q,0.05}}{\sqrt{1+n_0/n}}
	\right). \]

Of course, for any fixed $n$, we find again the same predictive
probabilities. This is due to the following fundamental property of
the cdf (Lecoutre, 1999):

\[ \Pr\Big(K'_{q_0,q}(a)<x\Big) = \Pr\Big(K'_{q,q_0}(x)>a\Big). \]

The K-prime distribution can also be used to make predictive 
statements about the standardized difference $d/s$ in a future sample.
In the same situation as above (two groups with a same sample size) 
we have:

\[ \frac{d}{s} \sim \sqrt{\frac{2(n_0+n)}{n_0\ n}}\ K'_{q_0,q}
\left(\frac{d_0}{s_0}\sqrt{\frac{n_0\ n}{2(n_0+n)}}\right). \]

When $q \rightarrow \infty$, this distribution tends to the distribution 
of the parameter $\delta/\sigma$. Thus, with a very large value of $n$, 
it could be used to get a statement about the \textit{population} 
standardized difference (as an alternative to the $\Lambda$-prime cdf).

For instance, suppose that $d_0/s_0=3$ and $n_0=100$. Then, taking $n=500000$, 
\[ Pr(d/s>2.731804)=1-Pr\left[ K'_{198,999998}
\left(21.21108\right)<19.31484\right]=0.9000. \]

But, actually, KPRIMECDF cannot provide a sufficiently accurate answer, 
even when the maximum absolute error parameter is set to $10^{-2}$, 
and issues an error message, while the algorithm without round-off 
error calculation returns a value (0.92) which is in error by $2$ times 
the required accuracy.

Concerning the K-square distribution, it can be used for the sample
size determination in ANOVA designs. For instance, a simple $g$-sample
experiment is designed to test the equality of $g$ means. A pilot
study has already been conducted with $g$ groups of equal sample size
$n_0$, and a $F$ ratio $F_0$ has been obtained (under the usual
normal model). Assuming an initial non informative prior, the
posterior predictive distribution for the $F$ ratio in the planned
experiment with $n$ subjects in each group is a K-square distribution:

\[ F \sim \frac{1+n/n_0}{g-1}\ K^2_{g-1,gn_0-g,gn-g}
\left(\frac{g-1}{1+n_0/n}F_0\right).\]

Suppose that $g=3$, $n_0=10$, $F_0=3.6$ and $n=30$. Then, given the
first results, $F$ is distributed as $2K^2_{2,27,87}(5.4)$ and the
probability of obtaining a significant $F$ test at 0.05 level is
$\Pr(F>3.1013) = 0.7792$. In order to get predictive probabilities
equal to 0.80 and to 0.90, $n=33$ and $n=54$ subjects in each group
are respectively needed.

Other uses of the K-prime and K-square distributions are the
computation of the cdf of the sampling distributions of correlation
coefficients. The cdf of the sample coefficient $r$, involving a
sample of $n$ independent observations from a bivariate normal
population with population coefficient $\rho$, is a particular case
of the K-prime distribution:

\[ \Pr(r<x) = \Pr\left[K'_{n-1,n-2}
\left(\sqrt{n-1}\,\frac{\rho}{\sqrt{1-\rho^2}}\right) <
\sqrt{n-2}\,\frac{x}{\sqrt{1-x^2}} \right]. \]

\noindent
The cdf of the square of the sample coefficient $R^2$, involving a
sample of $n$ independent observations from a $p$-variate normal
population with square multiple correlation coefficient $\rho^2$,
is a particular case of the K-square distribution:

\[ \Pr(R^2<x) = \Pr\left[K^2_{p-1,n-1,n-p}
\left((n-1)\frac{\rho^2}{1-\rho^2}\right) <
\frac{n-p}{p-1}\ \frac{x}{1-x}\right]. \]

\section{Limitations and possible improvements} \label{secPBI}

Drawbacks of Method~1 algorithms \citep[in the terms of ][]{BK03} 
led to the development of Method~2 algorithms. In Method~1, the 
iterations start at index $j=0$ which maximizes $H_j(x)$, while in 
Method~2 they start at index $j=k$ which maximizes $g_j$. 
Nevertheless, the latter is not systematically better. For instance, 
it can happen that the initial recurrence increment for the $H_j$'s 
is too small with respect to the machine limit so that a zero is 
returned and recurrence is impossible: e.g., 
for the K-square cdf, this increment term is lower than $10^{-307}$ 
when $p=10, q=20, r=30, a^2=500$ and $x=0.1$. More generally, whenever
$H_k(x)$ tends to zero quickly with respect to $k$, Method~1 algorithms
perform better than Method~2 algorithms, because only the first terms
of the series (\ref{EQ1}) contribute significantly to the sum. And 
when $H_k(x)$ is still close to $H_0(x)$, Method~2 is quasi optimum 
(with the same parameters as in the preceding example, this is the 
case when $x=99: H_0(99)\approx1$ and $H_{225}(99)=0.994$).

Obviously, the best method would be to start iterations at the index
(between $0$ and $k$) which maximizes the product $g_jH_j(x)$ and 
not only one of the terms. However, this is not easy to
determine in general. A tempting solution, when $H_k(x)$ is considered
too small, would be to choose the
modified index, say $k'$, such that $H_{k'}(x)$ reaches a
predetermined value; unfortunately, such an inversion of the beta cdf
involves an iterative procedure and so is to be discarded on grounds
of speed efficiency. As an alternative, we propose to simply
lower $k$ by multiplying it by the argument of the incomplete beta
function ($px/(px+r)$ for the K-square and $x^2/(x^2+r)$ for the
K-prime).

For example, for the distribution $K^2_{10,80,200}(500)$, when $x$
takes the values 35, 30, 25, and 22, the number of iterations is 
always 202 (for a precision of $10^{-4}$), while when turning to 
the modified starting index, it drops respectively to 155, 146, 
136 and 128.

\section{Concluding remarks} \label{secCLR}

We presented an algorithm for two Bayesian predictive distributions
of importance for monitoring experiments. This algorithm
includes round-off error calculation and is applicable to any
cumulative distribution function that can be expressed as a discrete
mixture of continuous distributions such that the recurrence relation
for the discrete coefficients is multiplicative and the recurrence
relation for the continuous distribution is additive. However, this 
kind of error calculation (which is only an approximation, of course) 
is time consuming, and when speed is a crucial factor, 
it has to be introduced only when deemed necessary. 
It will be the case, for example, when the required accuracy is high 
and/or when the number of iterations is large so that the precision 
of the result may be suspected. 
In this regard, the material used (computer and compiler) is of 
importance, particularly through the variable noted $\epsilon$, 
the precision of an ``elementary'' recurrence calculation. 
For instance, two different computers/compilers storing variables into
words of the same size could have different $\epsilon$ if they use
registers of different size to perform computations.
We also considered the case where the starting index of iterations is
such that recurrence is impossible due to underflows. The proposed 
solution, which is an approach to the problem of finding the 
optimum starting index, is to lower this index by a quantity which is
the argument of the incomplete beta function, a choice we made on
empirical grounds and that is likely to be improved.



\begin{thebibliography}{}


\bibitem[Benton and Krishnamoorthy(2003)]{BK03}
Benton, D., Krishnamoorthy, K., 2003. Computing discrete mixtures of
continuous distributions: noncentral chisquare, noncentral $t$ and
the distribution of the square of the sample multiple correlation
coefficient. Comput. Stat. Data An. 43, 249-267.

\bibitem[Frick(1990)]{Frick}
Frick, H., 1990. A remark on Algorithm AS 226: Computing noncentral
beta probabilities. Appl. Statist. 36, 311-312.

\bibitem[Johnson, Kotz and Kemp(1993)]{JKK93}
Johnson, N.L., Kotz, S., Kemp, A.W., 1993. Univariate Discrete
Distributions. Second edition. Wiley, New York.

\bibitem[Lecoutre(1984)]{BL84}
Lecoutre, B., 1984. L'Analyse Bayésienne des Comparaisons. Presses
Universitaires de Lille, Lille.

\bibitem[Lecoutre(1999)]{BL99}
Lecoutre, B., 1999. Two usefull distributions for Bayesian predictive
procedures under normal models. J. Statist. Plann. Inference 79, 93-105.

\bibitem[Lecoutre(2001)]{BL01}
Lecoutre, B., 2001. Bayesian predictive procedure for designing and
monitoring experiments. In Bayesian Methods with Applications to
Science, Policy and Official Statistics, Luxembourg: Office for
Official Publications of the European Communities, 301-310.

\bibitem[Lecoutre, Guigues and Poitevineau(1992)]{psi2}
Lecoutre, B., Guigues, J.-L., Poitevineau, J., 1992. Distribution of
quadratic forms of multivariate Student variables. Appl. Statist. 41,
617-627.

\bibitem[Majumder and Bhattacharjee(1973)]{AS63}
Majumder, K.L., Bhattacharjee, G.P., 1973. The incomplete beta 
integral. Appl. Statist. 22, 409-411.

\end{thebibliography}
\end{document}